\documentclass[12pt]{amsart}
\usepackage{amsmath,amssymb,url}

\newcommand\Aut{\mbox{\rm Aut}}

\newcommand\gen[1]{\langle#1\rangle}
\newcommand\sz[1]{\left|#1\right|}

\def\GAP{{\sf GAP}}

\def\A{\mathcal{A}}
\def\B{\mathcal{B}}
\def\R{\mathcal{R}}
\def\Ess{\mathcal{S}}
\def\T{\mathcal{T}}
\def\V{\mathcal{V}}
\def\X{\mathcal{X}}
\def\Y{\mathcal{Y}}

\author{Alexander Hulpke}
\address{Department of Mathematics, Colorado State University, 1874 Campus
Delivery, Fort Collins, CO, 80523-1874}
\email{hulpke@colostate.edu}

\dedicatory{In memory of Richard Parker, who taught us how sophisticated
computer programs can tame the largest of groups}

\title{Arithmetic in Group Extensions}

\begin{document}
\maketitle

\begin{abstract}
We describe a generalization of the concept of a Pc presentation that
applies to groups with a nontrivial solvable radical. Such a representation
can be much more efficient in terms of memory use and even of arithmetic,
than permutation and matrix representations. We illustrate its use
by constructing a maximal subgroup of the sporadic Monster
group and calculating its -- hitherto unknown -- character table.
\end{abstract}

PcGroups~\cite{SOGOS},
which are groups given by a polycyclic presentation with elements
represented as words in generators in normal form, have been one of the
success stories of Computational Group Theory.

The reasons for this success are multiple:
\begin{enumerate}
\item The representation applies to an important class of groups -- finite
solvable groups -- and many of them do not have
faithful permutation or matrix representations of small degree that could be
used alternatively.
\item Group elements can be represented effectively
(elements can be represented by exponent vectors in at most
$2\log_2(\sz{G})$ bits per element, and there are about
$\log(\sz{G})^2$ many polycyclic relations).
The {\em collection} process of rewriting elements into normal form provides an
effective multiplication routine~\cite{NewmanNiemeyer15}.
\item Reduction to factor groups $G/N$ for $N$ elementary abelian provides
an inductive paradigm for effective algorithms for many tasks
(e.g.~\cite{cellerneubueserwright,SOGOS,meckyneubueser})
in which
calculations are reduced to linear algebra and orbit/stabilizer calculations.
\item Many quotient algorithms (such as the $p$-Quotient~\cite{pquotient} or
Solvable Quotient~\cite{niemeyer94,plesken87}), which can find certain
quotients of finitely presented groups, internally use properties and
features of PcGroups for efficiency; they also produce output in the form of
PcGroups.

Indeed, the construction of groups of small order~\cite{bescheeickobrien02}
can be, broadly,
considered as an application of such quotient algorithms, namely to
a free group.
\end{enumerate}
These properties are not independent; for example, if the arithmetic was
not effective, there would be little point in using quotient algorithms, and
without quotient algorithms it would be much harder to come up with examples of
groups in polycyclic presentations.

This success motivates the desire to generalize these concepts to a wider
class of groups, ultimately reaching the class of polycyclic-by-finite
groups, which is natural to consider for algorithmic
questions~\cite{baumslagcannonitorobinsonsegal}.
\medskip

In such a generalization, we consider a group $G$ with a solvable normal
subgroup $B\lhd G$ as an extension
of a finite factor group $G/B$ by the normal subgroup $B\lhd G$. We
represent $B$ as PcGroup with a polycyclic generating set, and $G/B$
with a faithful permutation
representation. We shall call such a combination
a ``hybrid'' representation. While there is no formal requirement for
$B$ to be finite, we shall only study the finite case in this paper to avoid
dealing with issues of large integer exponents.
\bigskip

Indeed, such generalizations of concepts for PcGroups
have been pursued for some time. The {\em
Solvable Radical} paradigm~\cite{coxcannonholtlatt} uses exactly such a
structure to produce algorithms for permutation and matrix
groups, often generalizing earlier work for solvable groups.
Other examples are given by the construction of lists of perfect
groups~\cite{holtplesken89,hulpkeperfect}, which
produce the groups as
extensions of a known factor group by a solvable normal subgroup.

A hybrid representation of elements (and associated arithmetic) for finite
groups has been studied before, e.g. in~\cite{sinananholt}.
\medskip

Beyond the theoretical interest in representing a larger class of groups in
a uniform way, the highly memory-efficient element storage makes such a
representation interesting for working with groups for which the only
faithful permutation or matrix representations are of large degree: A
permutation of degree $10^6$ requires 3.8MB of storage, making permutation
representations of degree several million quickly infeasible for practical
use,
even on larger computers. This holds, for example, for groups constructed in
the first place as
extensions, which (when not solvable) have hitherto often been inaccessible
for computations. An example is the group constructed later in
Section~\ref{monmax}.
\smallskip

Another potential application of such a representation is in computing
automorphism groups~\cite{holtcannonautgroup}: The standard algorithm
computes the automorphism group of a group $G$ from the automorphisms of a
factor $G/N$ for $N\lhd G$ elementary abelian. The hardest part of such a
calculation
is often to determine the subgroup $S\le \Aut(G/N)$ that consists of
the automorphisms of $G/N$ that are induced by
automorphisms of $G$. Since this is a property that can be tested for an
individual automorphism, the construction of $S$ is similar to a stabilizer
calculation, and thus requires a means of determining the order of
subgroups of $\Aut(G/N)$, for which no low-degree permutation representation
might be known (or even exist).
On the other hand, the way $\Aut(G/N)$ will have been computed,
immediately produces a large normal
subgroup (the group $C$ in the labeling of~\cite{holtcannonautgroup}), which
is solvable, making such a hybrid representation potentially attractive.
\medskip

Concerning quotient algorithms, \cite{dietrichhulpke21} describes a general
``hybrid'' quotient algorithm that lifts a known finite quotient to an
extension with a module. The implementation of this algorithm, and the
construction of extensions used therein, have motivated the work described
in this article.

We shall describe such ``hybrid'' representations of groups as formal
extensions in which arithmetic is sufficiently fast for practical
calculations.
The approach has been implemented in the system~{\GAP} \cite{GAP4}
and can make use, with minimal changes, of existing implementations of {\em
Solvable Radical}-style algorithms.

It is now utilized in the implementation of the algorithm
from~\cite{dietrichhulpke21}, and we shall describe its use in determining the
hitherto unknown character table of a maximal subgroup of structure type
$2^{10+16}.O_{10}(2)$ in the sporadic simple Monster group,

\section{Hybrid Group Arithmetic}

A group extension is not determined solely by the isomorphism types of the
normal
subgroup and the factor group. A convenient way of describing the extra "glue" in
a concise way is through the concept of (confluent) rewriting
systems~\cite{knuthbendix,simsbook}. In particular, this allows us to
consider the new representation as a generalization of Pc presentations
(which are a special kind of confluent rewriting system).

We therefore start by summarizing some language:
For a generating set $\X$ of a group, denote by $w(\X)$ a word in this
generating set. A {\em rewriting system} over
$\X$ is a set $\R$ of pairs of words ({\em rules}) in $\X$, each rule
$(l,r)\in\R$
considered to consist of a {\em left} part $l$ and a {\em right} part $r$.
(We shall assume that $r$ is smaller than $l$ in a suitable ordering; this
is required to ensure the termination of reduction processes.)
We often write $l\to r$ to indicate that rules are intended to be used to
{\em reduce} a word $w(\X)$ by replacing an occurrence of a substring $l$
with a substring $r$. Such a process can be iterated for a set of rules,
resulting in a {\em reduced form} word in which no left part of any
rule occurs as a substring.
A rewriting system is called {\em confluent}, if every word has a unique
reduced form, that does not depend on the order in which (and positions in
the word at which) rules are applied.

A rewriting system can be considered as yielding a monoid presentation by
interpreting a rule $l\to r$ as identity $l=r$. If the system is confluent,
elements in normal form correspond exactly to elements of the presented
monoid.

Since a group is a monoid, it makes sense to consider rewriting systems that
give presentations for particular groups. While this generally requires
considering
generator inverses as new generators,
in the case of finite order generators, this
can be achieved by adding rules $x^{o(x)}\to 1$ for each generator $x$,
where $o(x)$ denotes the order of the element represented by $x$.

\subsection*{The constituent parts of a hybrid group}

We define a {\em hybrid group} $G$ as an extension of a polycyclic normal
subgroup $B\lhd G$ by a finite factor group $A=G/B$, with natural
homomorphism $\nu\colon G\to A$ and generators and
relations that reflect the extension structure:

Similar to \cite[\S 7]{dietrichhulpke21} the generators for $G$ consist of two
sets $\A=\{a_1,\ldots,a_k\}$ and $\B=\{b_1,\ldots b_l\}$ with $\A,\B\subset
G$. While there is no formal restriction on $\A$, we require that $\B$ is
a polycyclic generating set for $B$.
We also introduce two abstract
generating sets $\X=\{x_1,\ldots x_k\}$ corresponding to $\A$,
and $\Y=\{y_1,\ldots,y_l\}$ corresponding to $\B$.
We further assume a faithful permutation representation $\varphi\colon A\to
S_\Omega$,
which shall be given in the form of permutation images of the generators in $\A$.

We denote the polycyclic relations for $B$, considered as rewriting rules in the
formal generator set $\Y$, by
$\R_1$, and we thus talk about words in $\Y$ in normal form.
(In practice, we represent the group $B$ as an actual PcGroup on the computer.)
\smallskip

Next, we assume that we have a confluent rewriting system
$\R_A$ in the generators $\nu(\A)$, with
rules written in the formal symbols $\X$. Such a rewriting system can be
obtained, e.g. as described in~\cite{schmidt10}, and typically we will
choose the generating set $\A$ based on this rewriting system.

We modify these rewriting rules so that they reflect the extension structure
for $G$; this produces rules analogous to the power relations of a Pc
presentation:

Consider a rule in $\R_A$ of the form $l_i(\X)\to r_i(\X)$.
Evaluating the words $l_i(\X)$ and $r_i(\X)$ in the generators $\A$ yields
two elements of $G$ that have the same image in $A$ (since these are rules
for $A$). These two elements thus will differ by an element
$m_i=r_i(\A)^{-1}\cdot l_i(\A)\in B$, and this element $m_i$ can be written
as a normal form word $m_i=w(\B)$ in the polycyclic generating set. We
can thus write down a new rule
$l_i(\X)\to r_i(\X)\cdot w(\Y)$ in the abstract generators that describes
the behavior in $G$. The set of these
new relations is denoted
by $\R_3$. (In practice, we keep the relation set $\R_A$ as it was given and simply store for
each relation the corresponding cofactor $w(\Y)$ as a word in normal form.)
\smallskip

We finally have, for each generator $a_i$ of $A$ an automorphism
$\alpha_i\colon B\to B$ that describes the conjugation action of the element
$a_i$ on $B$, that is $b^{a_i}=\alpha_i(b)$.
This gives relations (analogous to the conjugation or commutator relations of a
Pc presentation) of the form
$y_j^{x_i}\to w(\Y)$ whenever $b_j^{a_i}=\alpha_i(b_j)=w(\B)$
in $G$, with the word in normal form for $B$. The set
$\R_2$ collects these relations.
\medskip

Using the same argument as in \cite[\S 7]{dietrichhulpke21}, one shows that
$\R_G:=\R_1\cup\R_2\cup\R_3$ is a confluent rewriting system for $G$ with respect to
a wreath product ordering.
\smallskip

The multiplication in $G$ then could be done by forming the product of words
in the free group and converting the result into normal forms with respect
to this rewriting system.
\medskip

In practice, for reasons of performance,
we handle the three classes of rewriting rules in different
ways: Rewriting rules in $\R_1$ are applied implicitly through the internal
representation of PcGroups.
This will typically be a Pc group collection process, but could equally well
use other methods, such as~\cite{leedhamgreensoicher98}.
A few extra remarks are in order regarding how we apply the rules
in $\R_2$ and $\R_3$:

In general, we follow a standard
``from-the-left`` strategy \cite[\S 2.4]{simsbook},
\cite{leedhamgreensoicher90}:

We can assume that the expression to be brought into normal form is a sequence
of letters from $\X$ and elements of $B$ (as PcGroup elements), the latter
in normal form. This will be ensured by a policy that whenever we encounter two
elements of $B$
in sequence, we replace them immediately by their product.
Doing so prioritizes the Pc rules and utilizes the
effective implementation of PcGroups.
\smallskip

For simplicity, we shall assume that the set $\R_A$ of rules for $A$ is
reduced, so that at any position in a word there is at most one rule that
could be applied, starting at the given position. (The same thus also holds
for the rules in $\R_3$.)

At any particular position of a word in $\X$, we can thus identify the
unique rule (if any) that would apply at this position by storing the left
sides of rules in the form of a prefix tree. Doing so avoids a costly
search through substrings.
\smallskip

In the process of bringing a word into normal form, we initially ignore the
$B$-parts and only consider the letters from $\X$. In this projection, we
then
determine the first
position at which the left side of a rule from $\R_3$ applies to the word and
retract that section to the original word expression.

Assume that this piece of the expression is
\[
E=
x_{j_0}\prod_{i=1}^{m} d_i x_{j_i}
=x_{j_0}d_1x_{j_1}\prod_{i=2}^{m} d_i x_{j_i}
\]
with $x_i\in\X$, $d_i\in B$ and $ls:=\prod_{j=0}^m x_{j_i}$ the left side
of a rule in $\R_3$.
We now iteratively move the $B$-elements to the right (while keeping the
value of the expression equal), so that the expression will contain
$ls$ verbatim and we can apply the corresponding rule:

First, we note that (from rules in $\R_2$, applied through the automorphism
$\alpha_{j_1}$ that represents conjugation by $x_{j_1}$)
we have that $d_1x_{j_1}=x_{j_1}d_1^{x_{j_1}}=x_{j_1}\alpha_{j_1}(d_1)$.
This results in the word
\[
x_{j_0}x_{j_1}e_1x_{j_2}\prod_{i=3}^{m} d_i x_{j_i}
\]
with $e_1=\alpha_{j_1}(d_1)\cdot d_2$, immediately evaluated in normal form in
$B$. This expression has one fewer $B$-term in front of a letter from $\X$
while its $\X$-projection remains unchanged.

We repeat the process by moving $e_1$ past $x_{j_2}$ in the same way, yielding
\[
x_{j_0}x_{j_1}\underbrace{\alpha_{j_2}(e_1)\cdot
d_3}_{=:e_2}x_{j_3}\prod_{i=4}^{m} d_i x_{j_i}.
\]

(Should $d_2=1$ have been the identity, we could immediately have
applied the product $\alpha_{j_1}\cdot\alpha_{j_2}$, skipping two letters.)

Iterating the process, we eventually obtain
$(\prod_i x_{j_i})\cdot e=ls\cdot e$ with
$e\in B$ where $ls$ is the left
side of a rule $l_p(\X)\to r_p(\X)\cdot m_p$ in $\R_3$. We thus replace
the expression $E$ by $r_p(\X)\cdot(m_pe)$, evaluating the product
$m_pe$ in $B$.
\smallskip

This process repeats until the word has been brought into a form $w(\X)\cdot b$
with $w(\X)$ and $b$ both in normal form. This is a normal form with respect
to $\R_G$.
\medskip

The repeated evaluation of automorphisms is the most time-critical part of
this routine. To improve performance, we can split $\B$ into segments
(preferably aligned to a chief series) and enumerate all normal form words
that can be formed with generators from only one segment. Because of the
compact storage of polycyclic elements, it is easy to store all these products
and to cache, for every automorphism involved, the images of all segment
products.
Evaluating the image of an element $b$ under such an automorphism splits $b$
into segments and uses the cached images for the segments.
An elementary abelian layer on the bottom can be
handled effectively through representing the action by matrices, acting on
exponent vectors.
\smallskip

We also note that in some cases, multiple letters from $\X$ can be skipped
simultaneously by
applying a product of automorphisms.  This can be sped
up by caching for short (say $s\le 4$) word expressions
$\prod _{i=1}^s a_{j_i}$ the corresponding automorphism products
$\prod_{i=1}^s \alpha_{j_i}$, and then applying these products if possible.
\medskip

The inverse of a word $a\cdot b$ with $a=\prod a_{j_i}$ and
$b\in B$ will be $b^{-1}\cdot a^{-1}$. To determine $a^{-1}$ we form the inverse word for $a$ (reverting ordering and
inverting letters)
and bring this word into normal form $\prod a_{k_i}\cdot t$. This shows that
$t^{-1}\cdot \prod a_{k_i}$ will be the inverse of $a$ and
$b^{-1}\cdot t^{-1}\cdot  \prod a_{k_i}$ the inverse of the whole word.

Since this process requires collection, and since the formally inverse word
typically is far from normal form, it is, in
fact, cheaper to precompute (with exactly the described method) the inverses for
every generator $a_i\in \A$, and to form the inverse of $a$ as a product of these
inverses, evaluating this longer product in a single collection process.
Again, one can trade space for speed by caching inverses of short generator
products.
\bigskip

To determine the order of an element $g\in G$, we first find the order $o$ of the
permutation image $\nu(g)$ and then the order of $g^o$ in the polycyclic
group $B$.

\bigskip
In comparison with the setup of~\cite{sinananholt}, a main difference of
this approach is in
requiring a confluent rewriting system for the factor group. In the
applications we primarily considered, this is not an obstacle, since the construction
of the group $G$ in the first place (as a formal extension, or through a
quotient algorithm) already required such a rewriting system.
\medskip

Note that we are not requiring $B$ to be as large as possible, nor do we
require
$A$ to have a trivial Fitting subgroup. But if this is not the case, we can
compute the solvable radical $S\lhd A$ using the permutation representation
and from this build a polycyclic presentation for the subgroup $R\lhd G$ that is
the full pre-image of $S$ in $G$. We thus get data structure information for
an isomorphic hybrid group with a normal subgroup isomorphic to the radical $R$
and factor group $A/S$.
The arithmetic in this new hybrid representation is
typically faster than in the original group, and it can be
advantageous, as far as arithmetic performance is concerned, to perform
arithmetic through the use of such a ``shadow''
representation.
\medskip

Homomorphisms from $G$, given by images of the standard generators, are
evaluated easily. This holds in particular for the natural homomorphism
$\nu\colon G\to G/B\cong A$ with the image in a permutation representation.
Other generating sets are handled through the more general approach for
subgroups we shall describe next.

\section{Subgroups, Factor Groups, and Homomorphisms}

The only extra feature needed to allow us to use many existing
group-theoretic algorithms is a membership test in
subgroups, as well as compatible data structures that allow for calculations
in such subgroups. This compatible data structure will also provide
information that will let us easily determine the subgroup order.

We shall describe how to do this for a subgroup $S\le G$,
given by a set $\Ess$ of generators.

We first calculate generator images $\nu(\Ess)$ under the natural
homomorphism and, using standard permutation group methods, determine a
presentation in these images.
(Standard techniques, such as straight-line programs, can be used to avoid
problems caused by long absolute word lengths.)
The relators of this presentation, evaluated in $\Ess$, form a set $L_0$
that generates a subgroup $\gen{L_0}\le S$,
whose normal closure in $S$ is $S\cap B$.
We then extend $L_0$ to a set $L$ of subgroup generators of $S\cap B$
by initializing $L:=L_0$ and
systematically forming $S$-conjugates (conjugating every element of $L$ by
every generator in $\Ess$), adding new elements to $L$ if they are not yet in
$\gen{L}$,
and iterating until no new conjugates arise.

The required membership tests in $\gen{L}$ are performed in the polycyclic
group $B$ and thus are cheap. Indeed, we can have this process create an
induced generating set (IGS,~\cite{SOGOS}) for $S\cap B$, as well as word
expressions in $\Ess$ for every element of this IGS.

For performance reasons, it can be advantageous to first search
systematically for short word expressions that lie in the kernel, typically
resulting in overall shorter words.

Together, the generating set $\Ess$, the
permutation group structure for $\gen{\nu(\Ess)}=\nu(S)$,
and the word expressions
for the IGS for $S\cap B$ are called the {\em hybrid bits} of $S=\gen{\Ess}$
and relate the subgroup to the hybrid structure of its parent group.

Expressing the original generators $L_0$ in terms of the IGS provides a way
to determine a presentation of $S$ in terms of $\Ess$.
\smallskip

Note that such an ``induced'' structure for $S\le G$
is not necessarily equivalent to a
hybrid structure for $S$, determined from scratch. This is because
$B\cap S$ might be strictly smaller than the radical of $S$. We also only
construct a presentation, and not a confluent rewriting system for
$\nu(S)$, though this could easily be done if it was desired to represent
$S$ (on its own) as a hybrid group.
\smallskip

To test membership of $g\in G$ in $S$, we first test membership of
$\nu(g)$ in $\nu(S)$. If this necessary condition holds, the standard
permutation membership test also gives a word expression of $\nu(g)$ in
$\nu(\Ess)$. Evaluating this word in $\Ess$, with the result $h\in G$, thus
reduces the membership test to determining whether $h^{-1}g\in S\cap B$,
which is done using the IGS for this subgroup.
\medskip

Implicitly, this process also
expresses $g\in S$ as a word in $\Ess$. Thus, the same process also allows us
to evaluate arbitrary homomorphisms, given as images of the generators $\Ess$.
\medskip

For two subgroups $U\le S\le G$, we determine a transversal for $U$ in $S$
as follows:
Let $\bar \T$ be a right
transversal of $\nu(U)$ in $\nu(S)$ and $\T=\{T_1,\ldots,T_p\}$
a corresponding list of pre-images in $S$.
(If $S=G$ these pre-images are immediately constructed. Otherwise,
we need to decompose into generator images in $\nu(S)$.)
We also form a transversal $\V=\{V_1,\ldots, V_q\}$ for $U\cap B$ in $S\cap B$.
Then the products $V_i\cdot \bar T_j$ form a transversal of $U$ in $S$:
Given $s\in S$, we determine the index $j$ such that $\nu(U)\cdot
\nu(s)=\nu(U)\cdot \bar T_j$.

Then $s\cdot T_j^{-1}\in (S\cap B)\cdot U$. But the transversal
$\V$ for $U\cap B$ in $S\cap B$ is also a transversal for $U$ in $(S\cap
B)U$, so we can find an index $i$ such that $s\cdot T_j^{-1}\cdot
V_i^{-1}\in U$. This both describes the transversal and provides a method for
identifying the coset in which an element of $U$ lies.

The quick identification of cosets in a transversal then is a crucial step
towards computation of permutation representations, allowing to turn
(subgroups of) hybrid groups into permutation groups.
\medskip

Factor groups $G/N$ for $N\lhd G$ can be represented as hybrid groups on
their own, using a similar approach. If $N\le B$, we
consider $G/N$ as an extension of $B/N$ by $G/N$, changing tails of the
rewriting rules from $B$ to $B/N$.
(In the special case $N=B$, we of course can simply take $\nu(G)$.)

More generally, $G/N$ is an extension of $B/B\cap N$ by $\nu(G)/\nu(N)$.
We add
word representations for generators of $\nu(N)$ to the
rewriting rules for $\nu(G)$, yielding rewriting rules for
$\nu(G)/\nu(N)$ which we make confluent again with the Knuth-Bendix
method~\cite{knuthbendix}. (This is typically much easier than starting a Knuth-Bendix process
from an arbitrary presentation.)
As for a permutation representation of $\nu(G)/\nu(N)$, we use the standard
heuristics in {\GAP} for determining such a representation. We note that the
typically $\nu(G)$ is a group with trivial radical. This means that its
structure is rather limited and lends itself to good
permutation representations.
\bigskip

Together, these methods encompass all representation-specific operations
that are required for many already existing algorithms that use the
{\em Solvable Radical} paradigm: These are algorithms for
conjugacy classes~\cite{cannonsouvignier}, (maximal)
subgroups~\cite{coxcannonholtlatt,holtcannonmaxgroup,eickhulpke01},
automorphism groups~\cite{holtcannonautgroup}, Sylow subgroups~\cite{eickhulpke12}, or normalizers~\cite{glasbyslattery}.
\medskip

A, perhaps surprising, place where it can be necessary to adapt existing
implementations concerns reductions in the number of generators.
Algorithms that originate in the world of permutation groups often just need
to store {\em some} generating set for a subgroup. Since storing
permutations is costly in memory, and since the cost of an orbit
calculation is proportional to the number of group generators, it is
convenient to keep the size of such a generating set small, for example by
removing redundant generators.
Algorithms that use Schreier generators thus often include tests to
eliminate redundant generators.

Such concerns do not hold as strongly for a (stabilizer) subgroup $S\le G$
of a hybrid group $G$ with radical $R$: Elements take
less memory, while an IGS for $S\cap R$ will in any case be part of the
data structure to store $S$.

In consequence, it is thus advisable to
revisit uses of subgroup generator reduction in implementations.

\section{Some Examples}

The author's implementation of such groups is currently available
on {\tt github}
under \url{https://github.com/hulpke/hybrid}. (The same file also provides
the functionality of~\cite{dietrichhulpke21}, since it now builds on hybrid
groups.)
The code currently requires the use of a development version of {\GAP} since
some algorithm implementations required minimal changes -- mostly cleaning up careless
representation-specific function calls.

It is the author's intention to make this code available in the future as part
of a standard {\GAP} distribution.
\medskip

In this implementation,
hybrid groups can be created
\begin{enumerate}
\item From an existing permutation or matrix group
\item As extensions, based on 2-cohomology information
\item As quotients of a larger group $H$, using a homomorphism $\varphi\colon
H\to P$ in a finite permutation group $P$, and a homomorphism
$\pi\colon\ker\varphi\to S$ into a solvable group.
\item As output of the hybrid quotient algorithm~\cite{dietrichhulpke21}.
\end{enumerate}
\smallskip

For these groups, a
significant number of already existing {\em Solvable Radical} paradigm
algorithms become available immediately (or with minimal changes).
\medskip

We compare this implementation of hybrid groups in {\GAP} with some of the
examples used in~\cite{sinananholt}. They are given in Table~\ref{timings}.

\begin{table}
\begin{tabular}{l|ccccc}
Group&$x\cdot y$&$x^{-1}$&$\mbox{Syl}_2$&$\mbox{Syl}_q$&Classes\\
\hline
$2^{10}\cdot M_{12}$&0.48&0.34&20&11: 8ms&68: 455ms\\
$2^{1+8}\cdot O^+_8(2)$&17&15&560&7: 152ms& 332: 105877ms\\
$2^{24}.A_5$&0.05&0.03&90&5: 13ms&236: 803ms\\
$13^{55}.A_5$&0.21&0.78&82&13: 12217ms&---\\
$7^{45}\cdot A_8$&0.25&0.08&1240&5: 30ms&---\\
$2^{10+16}.O^+_{10}(2)$&17&10.8& 4300&31: 510&478: 229818\\
\end{tabular}
\caption{Timings (in ms) for some hybrid group calculations}
\label{timings}
\end{table}

\begin{table}
\begin{tabular}{l|cccc}
Group&$x\cdot y$&$x^{-1}$&$\mbox{Syl}_2$&$\mbox{Syl}_q$\\
\hline
$2^{10}\cdot M_{12}$&7.1&8.35&12217&11: 8280ms\\
$2^{1+8}\cdot O^+_8(2)$&41&28&16680&7: 2ms\\
$2^{24}.A_5$&0.175&0.275&91&5: 71ms\\
$13^{55}.A_5$&0.22&0.4&156&13: 30ms\\
\end{tabular}
\caption{Comparison runtimes from~\cite{sinananholt}, scaled by factor 4}
\label{timingsSH}
\end{table}

Timings are on a Mac~Studio with an  M1~Max processor and 64GB memory
(Geekbench~6: 2400). We estimate, from the publication date and processor
specification of~\cite{sinananholt}, that our processor is roughly a factor
of 3-4 faster. We thus reproduce the runtimes from~\cite{sinananholt},
scaled by a factor 4, in Table~\ref{timingsSH}, which shows our
implementation to be competitive.

For comparison, multiplying permutations of degree $10^6$ takes about 1ms
each.

Sylow subgroup calculations and conjugacy classes
exclude the cost of building a {\em Solvable Radical} data structure
(which, once computed, will be used for all subsequent calculations).

We did not compare timings for the center, as we did not implement the
special method from~\cite[\S 5]{sinananholt}.

As for the groups in the examples, we note that the structure descriptions and
further information
in~\cite{sinananholt} do not necessarily identify a unique group up to
isomorphism: There are three pairwise non-isomorphic
extensions of nonsplit type $2^{10}{\cdot}M_{12}$, two of which have minimal
permutation degree 264. There also are multiple groups of structure
$2^{1+8}\cdot O^+_8(2)$.

In both cases, we picked one particular extension at random for the
experiments. (Comparing
a few run times for the other candidate groups fitting the same
description seems to indicate that a different choice of group would not
have grossly distorted the picture.)

The group $7^{45}\cdot A_8$ is a nonsplit extension with an irreducible
module of dimension 45, constructed as formal extension.

The last group, $2^{10+16}.O^+_{10}(2)$, is constructed in the following
section.

\section{The Character Table of a Maximal Subgroup of the Monster}
\label{monmax}

We consider the group $G=2^{10+16}.O^+_{10}(2)$, the 5th maximal subgroup of
the sporadic Monster group $M$ \cite{ATLAS}. It has order $2^{46}3^55^27\cdot
17\cdot31$. K.~Lux had suggested to the
author to construct this group from scratch, with the aim of obtaining its
(hitherto unknown)
character table.

While all maximal subgroups of $M$ have been constructed
recently~\cite{DietrichLeePisaniPopiel,DietrichLeePopiel},
these calculations are
subgroup generators in the specific format of the {\tt mmgroup}
package~\cite{mmgroup} and
as such  not necessarily amenable to calculations such as conjugacy classes.

Instead, we construct this group anew as a hybrid group: $O^+_{10}(2)$ has
two irreducible modules in dimension $16$ over $\mathbb{F}_2$, which are
swapped by an outer automorphism. The choice of either module thus yields
isomorphic extensions, and we therefore need to consider only one of the two
modules.
The associated 2-cohomology group has dimension $1$, thus there is
only one non-split extension to consider.
The required rewriting system for $O^+_{10}(2)$ was precomputed using the
approach from~\cite{schmidt10}, utilizing a BN-pair.

We construct this extension as a permutation group in a sequence of standard {\GAP} commands, as
described in Figure~\ref{figq}.
Suitable compatible permutation generators and corresponding matrix actions for
$O^+_{10}(2)$
are fetched from
the {\sf ATLAS} web pages~\cite{eATLAS}. The cohomology information yields a
presentation for the extension. A search for low index subgroups of a
subgroup mapping onto a point
stabilizer of $O^+_{10}(2)$ yields a permutation group
of degree 126480.
Its action on a suitable second maximal
subgroup reduces the permutation degree to 4590.
The whole construction takes about 6 hours on the
same machine as used in the previous section.

\begin{figure}
\begin{verbatim}
gap> g:=AtlasGroup("O10+(2)");; NrMovedPoints(g);
496
# Representation 9 is 16-dim matrices
gap> m:=AtlasGenerators("O10+(2)",9);;
gap> m:=GModuleByMats(m.generators,GF(2));;
gap> co:=TwoCohomologyGeneric(g,m);;
gap> co.cohomology; # check dim 1
[ <an immutable GF2 vector of length 21392> ]
gap> a:=FpGroupCocycle(co,co.cohomology[1],
> true);; # force perm rep
gap> p:=Image(IsomorphismPermGroup(a));;
gap> Size(p);
1540049859300556800
gap> NrMovedPoints(p);
126480
gap> m:=LowLayerSubgroups(p,2);; # second maximals
# avoid kernel
gap> m:=Filtered(m,x->not IsSubset(x,RadicalGroup(p)));;
gap> List(m,x->Index(p,x));
[ 1060991139840, 1300234240, 474300, 379440, 73440, 4590 ]
gap> rep:=FactorCosetAction(p,m[6]);; new:=Image(rep);;
gap> NrMovedPoints(new);
4590
\end{verbatim}
\caption{Constructing $2^{16}.O^+_{10}(2)$ in {\GAP}}
\label{figq}
\end{figure}

Calculating the 2-cohomology for this group and the unique 10-dimensional
module (the natural module of $O^+_{10}(2)$) yields a 2-dimensional
cohomology group. (This calculation took about 12 hours.)
We constructed the associated 4 extensions as
hybrid groups.

Using {\GAP}'s implementation of conjugacy classes~\cite{hulpkematclass},
we determine for each of these groups the conjugacy classes. This
calculation takes about 5 minutes for each group.
Table~\ref{classsz} gives some information that shows that the
groups must be pairwise nonisomorphic.

\begin{table}
\begin{tabular}{l|cc}
Candidate&Nr. Classes&Nr. Classes $|x|=4$\\
\hline
0 (split)&1063&289\\
1&478&62\\
2&1063&282\\
3&478&58\\
\end{tabular}
\caption{Distinguishing information of the constructed candidates}
\label{classsz}
\end{table}

We also observe that candidate groups 0 and 2 both have elements
of order $2$ with centralizer of order $2^{46}3^25\cdot 7\cdot 31$. Such a
group cannot embed into $M$, leaving us with two candidates (either of which,
based on
conjugacy classes, could potentially embed into $M$).
\medskip

It was thus decided to determine the character tables for the two remaining
candidates. This turned out to be somewhat challenging because of the large
order of the groups.
To speed up calculations, an attempt was made
to represent the groups as permutation groups. The best degree found
was $12,142,080$, obtained by acting on a conjugacy class of elements.
While we did not prove that this indeed is the smallest faithful degree,
a search
through iterated maximal subgroups indicates that this is most likely to be
the case. Attempts to work in such a permutation group failed, since the
storage needs for elements (a single permutation taking about 46MB of
memory) became prohibitive. We also note that the speed of element
arithmetic in the group (Table~\ref{timings}) is easily competitive with
permutations of degree $1.2\cdot 10^7$.
\bigskip

Instead, we proceeded via the hybrid representation, providing a challenge to
the implementation.
We implemented, experimentally, the
inductive approach from~\cite[\S 3]{ungerctbl} (that is, using induction and
lattice reduction, but not the modular approach of~\cite[\S 4]{ungerctbl})
in {\GAP}. Unfortunately, we found that in many examples this implementation
was unable to obtain irreducible characters.
We surmise that one reason behind this lies in the LLL implementation used.
While {\GAP} implements essentially a textbook version of the algorithm, the
version implemented in {\sf Magma} is understood to be highly optimized in
speed and performance.

\bigskip

We therefore took a classical character-theoretic approach to determine the character tables for the two remaining candidates for maximal subgroups. For each
candidate $G$, we first took the inflated characters of the factor
$2^{16}.O^+_{10}(2)$. We then computed a number of large subgroups (maximals
or maximals of maximals) $S\le G$ for which we were able to compute faithful
permutation representations, typically of degree $\sim 100,000$, from actions on
elements. We computed the character tables for these groups $S$ in this
permutation representation and induced
these characters to $G$. An LLL reduction of the resulting character pool
resulted in new irreducible characters. Replenishing the pool with tensor
products of the characters found so far and further lattice reduction then
ultimately
resulted in character tables for the two candidate groups.
These calculations took more memory than the author's desktop machine
provided and thus were done on a larger computing server that provided 1TB
of memory, at comparable processor speed. Calculations took
about 2-3 weeks for each group, in
a mix of automated calculations for maximal subgroups and manual work in
assembling the table for each $G$.
\bigskip

To determine the character tables of the large subgroups $S$, we
extended the author's {\GAP} implementation of the Dixon-Schneider
algorithm~\cite{dixonschneider,hulpkediplom} with an automatic version of
a similar character-theoretic approach: We start with character tables of
proper factor groups and inflate their irreducible characters. The
algorithm also produces, once it is known that the Dixon-Schneider algorithm
will require the calculation of class matrices for larger classes (more than
$10^5$ or $10^6$ elements), reducible
characters from tensoring of known irreducibles and induction from
subgroups; starting with cyclic subgroups and iterating to elementary
subgroups (in an echo of~\cite{ungerctbl}) and eventually maximal subgroups.
These reducible characters are again fed to an LLL reduction, in the hope of
obtaining
further irreducibles. For the (still very large)
groups arising here, this worked
significantly better than a pure Dixon-Schneider calculation and allowed us
to determine all desired character tables for a sufficient number of large
subgroups $S$ so that their induced characters eventually yielded the
character table of $G$.
\bigskip

Having obtained the character tables of the two candidate groups $G$,
we used the functionality of
{\GAP} and its character table library~\cite{GAPctbl}
to determine possible fusions from each table into the table of $M$.
(This functionality constructs putative fusions based on element orders,
centralizers, and power maps, with a subsequent verification whether induced
and restricted characters decompose into irreducibles with nonnegative
integer coefficients, which will eliminate many invalid options.)

Luckily, one table (candidate 3) proved to not allow for a fusion at all,
while the other candidate (number 1) allowed for exactly one fusion up to
table automorphisms. We thus found the desired character table by elimination.

This resulting character table has been included in the {\GAP} character table
library. A representation of it as a hybrid group can be obtained from
the author.
\smallskip

Added note in the revised version:
With this result, the only other maximal subgroup of the Monster whose
character table had not been known is of type
$2^{5+10+20}\cdot (S_3\times L_5(2))$. Since the submission of this paper,
this character table has been computed as well~\cite{pisanictbl},
using this same hybrid group representation. With this cited work, the
ordinary character tables of all maximal subgroups of the Monster are now known.

\section{Acknowledgements}
The author would like to thank the referees for careful reading and many
helpful remarks.

The calculations in section~\ref{monmax} were run on a large computing
server whose purchase in 2023 was supported by the Faculty Success program
of the College of Natural Sciences at CSU.

The author's work has been supported in part by
Simons Foundation Grant~852063,
which is gratefully acknowledged.

\bibliographystyle{amsplain}
\bibliography{mrabbrev,litprom}

\end{document}